\documentclass{article}
\usepackage{amsmath}
\usepackage{graphicx}
\usepackage{latexsym}
\usepackage{hyperref}
\usepackage{amsfonts}
\usepackage[all]{xy}
\begin{document}
\renewcommand{\theequation}{\arabic{section}.\arabic{equation}}
\newcommand{\be}{\begin{eqnarray}}
\newcommand{\en}{\end{eqnarray}}
\newcommand{\no}{\nonumber}
\newcommand{\ti}{\tilde}
\newcommand{\laa}{\lambda}
\newcommand{\la}{\langle}
\newcommand{\ra}{\rangle}
\newcommand{\ga}{\gamma}
\newcommand{\La}{\Lambda}
\newcommand{\ep}{\epsilon}
\newcommand{\de}{\delta}
\renewcommand{\thefootnote}{}
\newcommand{\pl}{\parallel}
\newcommand{\ov}{\overline}
\newcommand{\bet}{\beta}
\newcommand{\al}{\alpha}
\newcommand{\fr}{\frac}
\newcommand{\pa}{\partial}
\newcommand{\we}{\wedge}
\newcommand{\Om}{\Omega}
\newcommand{\na}{\nabla}
\newcommand{\D}{\Delta}
\newcommand{\ri}{\rightarrow}
\newcommand{\pp}{\phi_{\alpha i}}
\newcommand{\ii}{\int_{\Omega}}
\newcommand{\vs}{\vskip0.3cm}
\newcommand{\R}{I\!\!R^{n+1}}
\title {Eigenvalues of  the Wentzell-Laplace Operator and of the Fourth Order Steklov Problems \thanks{Partially supported by CNPq.}}
\footnotetext{
\hspace*{2ex}Key words and phrases: Eigenvalues, Wentzell-Laplace Operator, non-negative Ricci curvature, Steklov Problems, isoperimetric inequality, Euclidean ball.}
\author{ Changyu Xia and Qiaoling Wang } \date{}
\maketitle
\begin{abstract}
We prove a sharp upper bound and a lower bound for the first nonzero eigenvalue of the
Wentzell-Laplace operator on compact manifolds with  boundary and an isoperimetric inequality for the same eigenvalue in the case where the manifold is
a bounded domain in a Euclidean space.  We study some fourth order  Stekolv problems  and obtain isoperimetric upper bound for the first eigenvalue of them. We also find all the eigenvalues and eigenfunctions for two kind of  fourth order  Stekolv problems on a Euclidean ball.
\end{abstract}


\section{Introduction and Statement of the Results}
\markright{\hfill Eigenvalues of the Wentzell-Laplace Operator\hfill }
Let $n\geq 2$ and $(M, \la, \ra)$ be an $n$-dimensional compact Riemannian manifold with boundary.
We denote by $\overline{\Delta}$ and $\Delta$  the Laplace-Beltrami operators on $M$ and $\partial M$, respectively, and consider
the  eigenvalue problem for Wentzell boundary conditions
\be\left\{\begin{array}{l}
\ov{\Delta} u =0 \ \ \ \ \ \ \ \ \ \ \ \ \ \ \ \ \ \ {\rm in  \ }\ M,\\ -\beta\Delta u+\pa_{\nu} u=  \laa  u\ \ \ {\rm on  \ } \pa M,
\end{array}\right.
\en
where $\beta$ is a given real number and $\pa_{\nu}$ denotes the outward unit normal derivative. When $M$ is a bounded domain in a Euclidean space, the
above problem has been studied recently in [7]. A general derivation of Wentzell
boundary conditions can be found in [13].  Note that when $\beta=0$, the problem (1.1) becomes the
second order Steklov problem:
\be
\left\{\begin{array}{l}
\ov{\Delta} u =0 \ \ \  \ \ {\rm in \ \ } M,\\ \pa_{\nu} u=  p u \ \ {\rm on \ \ } \pa M,
\end{array}\right.
\en
which has a discrete spectrum consisting in a sequence
\be
\no
p_0=0<p_1\leq p_2\leq\cdots \rightarrow +\infty.
\en
When $\beta\geq 0$, the spectrum of the Laplacian with Wentzell condition consists in an
increasing countable sequence of eigenvalues
\be
\no\laa_{0, \beta}=0<\laa_{1, \beta}\leq\laa_{2, \beta}\leq\cdots\rightarrow +\infty,
\en
with corresponding real orthonormal (in $L^2(\pa M)$) eigenfunctions $u_0, u_1, u_2,\cdots.$
We adopt the
convention that each eigenvalue is repeated according to its multiplicity.  Consider the Hilbert space
\be
H(M)=\{ u\in H^1(M), {\rm  Tr}_{\pa M}(u)\in H^1(\pa M) \},
\en
where ${\rm Tr}_{\pa M}$ is the trace operator. We define on $H(M)$ the two bilinear forms
\be
A_{\beta}(u, v)=\int_{M}\ov{\nabla}u\cdot\ov{\nabla}v +\beta\int_{\pa M} \na u\cdot\na v, \ B(u, v)=\int_{\pa M}uv,
\en
where $\ov{\nabla}$ and  $\na$ are the gradient operators on $M$ and $\pa M$, respectively. Since  we assume  $\beta$  is nonnegative, the two bilinear forms are positive and the variational characterization for the $k$-th eigenvalue is
\be
\laa_{k, \beta}=\min\left\{\frac{A_{\beta}(u, u)}{B(u, u)}, u\in H(M), u\neq 0,\ \int_{\pa M} u u_i=0, i=0,\cdots,k-1.\right\}.
\en
When $k=1$, the minimum is taken over the functions orthogonal to the eigenfunctions
associated to $\laa_{0, \beta}= 0,$ i.e., constant functions. It is easy to see that if $\beta>0,\ p_1$ is the first non-zero eigenvalue of the Steklov problem (1.2) and $\eta_1$ the first non-zero eigenvalue of the Laplacian  acting on functions on $\pa M$, then we have
\be
\laa_{1,\beta}\geq \beta \eta_1+p_1,
\en
with equality holding if and only if any eigenfunction $f$ corresponding to $\laa_{1, \beta}$ is an eigenfunction corresponding to $p_1$ and
$f|_{\pa M} $ is an eigenfunction corresponding to $\eta_1$.

If $M$ is an $n$-dimensional Euclidean ball of radius $R$, then $\laa_{1, \beta}$ has multiplicity $n$, the corresponding  eigenfunctions are the coordinate functions $x_i,  i=1,\cdots,n,$  and
\be\laa_{1,\beta}=\fr{(n-1)\beta + R}{R^2}.\en
More generally all the eigenfunctions are the spherical harmonics and the eigenvalue associated to spherical harmonics of order $l$ is(Cf. [7])
\be\fr{k(k+n-2)\beta +kR}{R^2}.\en

In this paper we obtain a sharp upper bound for  the first nontrivial eigenvalue  of the problem (1.1).
\vskip0.2cm
{\bf Theorem 1.1.} {\it  Let $(M, \la, \ra)$ be an $n$-dimensional compact connected Riemannian manifold with boundary  and  $\beta$  a positive constant.  Denote by $\eta_1$ the first non-zero eigenvalue of the Laplacian  acting on functions on $\pa M$ and  $\laa_{1, \beta}$ the  first non-zero eigenvalue  of the problem (1.1). If the  Ricci curvature of $M$ is bounded below by $-\kappa_0$ for some non-negative constant $\kappa_0$ and the principle curvatures of $\pa M$ are bounded below by a positive constant $c$, then we have
\be
\laa_{1, \beta}\leq  \beta\eta_1+ \frac{2\eta_1+\kappa_0+\sqrt{(2\eta_1+\kappa_0)^2- 4(n-1)\eta_1 c^2}}{2(n-1)c}
 \en
 with equality holding if and only if $\kappa_0=0$ and  $M$ is isometric to an $n$-dimensional Euclidean ball of radius $\fr 1c$.}
 \vskip0.2cm
For bounded domains $\Om$ in a Euclidean space ${\mathbf R}^n$, we have an isoperimetric  upper bound for $\laa_{1, \beta}$ which depends only on the volume of $\Om$, the area of $\pa\Om$ and
the dimension $n$.
\vskip0.2cm
{\bf Theorem 1.2.} {\it  Let $\Om$ be a bounded domain with smooth boundary in ${\mathbf R}^n$ and $\beta$ a positive constant.   Then the  first non-zero eigenvalue $\laa_{1, \beta}$ of the problem
\be\left\{\begin{array}{l}
\ov{\Delta} u =0 \ \ \ \ \ \ \ \ \ \ \ \ \ \ \ \ \ \ {\rm in  \ }\ \Omega,\\ -\beta\Delta u+\pa_{\nu} u=  \laa  u\ \ \ {\rm on  \ } \pa \Om,
\end{array}\right.
\en
satisfies
\be
\laa_{1, \beta}\leq \fr{n|\Om|+ \beta(n-1)|\pa\Om|}{n|\Om|(|\Om|\omega_n^{-1})^{1/n}}
 \en
 with equality holding if and only if $M$ is isometric to a ball, where $|\Om|, |\pa\Om|$ and $\omega_n$ denote the volume of $\Om$, the area of
 $\pa \Om$ and the volume of a unit ball of ${\mathbf R}^n$, respectively.}
 \vskip0.2cm

Our next result provides a lower bound for the same eigenvalue.
\vskip0.2cm
{\bf Theorem 1.3.} {\it  Let $(M, \la, \ra)$ be an $n$-dimensional compact connected Riemannian manifold  and boundary $\pa M$ and $\beta$ a positive constant.   Assume that  the principal curvatures of $\pa M$ are bounded below by a positive constant $c$.

i) If the Ricci curvature of $M$ is bounded below by $-\kappa$ for some non-negative constant $\kappa$, then the first nonzero eigenvalue of the problem
(1.2) satisfies
\be
p_1>\frac{c\eta_1}{2\eta_1+\kappa},
\en
where $\eta_1$ is the first non-zero eigenvalue of the Laplacian  acting on functions on $\pa M$.

ii) If $M$ has non-negative Ricci curvature, then the first non-zero eigenvalue $\laa_{1, \beta}$ of the problem (1.1) satisfies
\be
\laa_{1, \beta} >\fr{\left(1+(n-1)c\beta+\sqrt{(n-1)^2c^2\beta^2+2(n-1)c\beta}\right)c}2
 \en
 }
\vskip0.2cm
In the case of $\kappa=0$, (1.12) has been proved by Escobar in [8]. An interesting question related to Theorem 1.3 is to find the best possible
lower bound for $\laa_{1, \beta}$ and $p_1$. Escobar conjectured that when $\kappa =0,\ p_1\geq c$ (Cf. [8]).  We believe that (1.13) could be improved as
\be
\laa_{1, \beta}\geq (n-1)\beta c^2+ c,
\en
with equality holding if and only if $M$ is isometric to a Euclidean ball of radius $\frac 1c$ in ${\mathbf R}^n$.

Now we come to eigenvalues of fourth order Steklov problems. Let $\Om$ be a bounded bounded domain in ${\mathbf R}^n$ and consider the problem
\be\left\{\begin{array}{l}
\ov{\Delta}^2 u- \tau \ov{\Delta} u =0 \ \ \ \ \ \ \ \ \ \ \ \ \ \ \ \ \ \ \ \ \ \  \ \ \ \  {\rm in \ \ } \Omega,\\ \fr{\pa ^2u}{\pa \nu^2}=0 \ \ \ \ \
\ \ \ \ \ \ \ \ \ \ \ \ \ \ \ \ \ \ \ \ \ \ \ \ \ \ \ \ \ \ \  {\rm on } \ \pa\Omega, \\
\tau \frac{\pa u}{\pa \nu}-{\rm div}_{\pa \Omega}(\ov{\nabla}^2u(\nu))-\frac{\pa \ov{\Delta}u}{\pa \nu}= \laa u \ \  {\rm on \ \ } \pa \Omega,
\end{array}\right.
\en
in the unknowns $u$ (the eigenfunction), $\laa$ (the eigenvalue), where $\tau > 0$ is a fixed positive
constant, $\nu$ denotes the outer unit normal to $\pa\Omega$, ${\rm div}_{\pa \Omega}$
 denotes the tangential divergence
operator, $\ov{\Delta}$ is the Laplacian on ${\mathbf R}^n$ and $\ov{\nabla}^2u$ the Hessian  of $u$. For $n= 2$, this problem is related to the study of
the vibrations of a thin elastic plate with a free frame and mass concentrated at the boundary.
The spectrum consists of a diverging sequence of eigenvalues of finite multiplicity
\be
0=\laa_0(\Om)<\laa_1(\Om)\leq\laa_2(\Om)\leq\cdots \leq\laa_j(\Om)\leq\cdots,
\en
where we  repeat the eigenvalues according to their multiplicity. The problem
(1.15) is the analogue for the biharmonic operator of the classical Steklov problem (1.2) for the Laplace
operator  and has been studied recently in [4]. The first nonzero eigenvalue of (1.15)  is usually called the fundamental tone and
can be characterized  by means of the Rayleigh principle
\be
\laa_1=\underset{\int_{\pa \Om}u=0}{\underset{0\neq u\in H^2(\Om)}{\min}}\frac{\int_{\Om} |\ov{\nabla}^2 u|^2+ \tau |\ov{\nabla} u|^2)}
{\int_{\pa \Om} u^2}.
\en
If $\Om$ is a ball of radius $R$ in ${\mathbf R}^n$, then $\laa_1=\fr {\tau}R$ (Cf. [4]).
The following  result provides  a sharp upper bound for the first nonzero eigenvalue of the problem (1.15) on bounded convex domains in ${\mathbf R}^n$.
\vskip0.3cm
{\bf Theorem 1.4.} {\it Let $\Om$ be a bounded smooth domain in ${\mathbf R}^n$ and assume that the principle curvatures of $\pa\Om$ are bounded below by a positive constant $c$. Let $\eta_1$ be the first  non-zero eigenvalue of the Laplacian  acting on functions on $\pa \Om$.
Then the first nonzero eigenvalue of the problem (1.15) satisfies
\be
\laa_1\leq \frac{\eta_1^2+ \tau\left(\eta_1+\sqrt{\eta_1^2-(n-1)\eta_1 c^2}\right)}{(n-1)c}-c\eta_1,
\en
with equality holding if and only if $\kappa=0$ and $\Om$ is a ball of radius $\frac 1c$.}
\vs
Another Steklov problem for the bi-harmonic operator we are interested is as follows:
\be
\left\{\begin{array}{l}
 \ov{\Delta}^2 u= 0  \ \ {\rm in \ \ } \Om, \\
\pa_{\nu}u = \pa_{\nu}(\ov{\D} u)+\xi u =0 \ \ {\rm on \ \ } \pa \Om.
\end{array}\right.
\en
This problem describes the deformation $u$ of the linear elastic supported plate $\Om$ under the
action of the transversal exterior force $f(x) = 0, \ x\in\Om $ with Neumann boundary condition $\pa_{\nu}u|_{\pa \Om}=0$(see,
[25] or p. 32 of [24]) and  was first discussed by J. R. Kuttler and V. G. Sigillito in [16].
The eigenvalues of (1.19) can be arranged as
\be
0=\xi_0<\xi_1\leq\xi_2\leq\cdots\ri +\infty.
\en
The first nonzero eigenvalue is given  by the following Rayleigh-Ritz formula.
\be
\xi_1=\underset{\int_{\pa \Om}u=0, \pa_{\nu}u|_{\pa\Om}=0}{\underset{0\neq u\in H^2(\Om)}{\min}}\frac{\int_{\Om} (\ov{\D} u)^2}
{\int_{\pa \Om} u^2}.
\en
When $\Om=B$(the unit ball in ${\mathbf R}^n$) we shall determine explicitly all the eigenvalues of (1.19). For each $k=0, 1,\cdots,$ let $\mathcal{D}_k$ be the space of harmonic homogeneous polynomials in ${\mathbf R}^n$ of degree $k$ and denote by $\mu_k$ the dimension of $\mathcal{D}_k$. We refer to [2] for the basic properties of  $\mathcal{D}_k$ and $\mu_k$. In particular, we have
\be\no
& &\mathcal{D}_0={\rm span}\{1\}, \ \ \ \mu_0=1,\\ \no
& &\mathcal{D}_1={\rm span}\{x_i, i=1,\cdots, n\}, \  \ \ \mu_1=n,\\ \no
& &\mathcal{D}_2={\rm span}\{x_ix_j, x_1^2-x_h^2,  1\leq i<j\leq n, h=2,\cdots, n\}, \ \ \ \ \mu_2=\fr{n^2+n-2}2.
\en
{\bf Theorem 1.5.} {\it If $n\geq 2$ and $\Om =B$, then we have

 i) the eigenvalues of (1.19) are $\xi_k=k^2(n+2k), k=0, 1, 2, \cdots $;

 ii) for all $k=0, 1, 2,\cdots,$ the multiplicity of $\xi_k$ is $\mu_k$;

 iii) for all $k=0, 1, 2,\cdots,$ and  all $\phi_k\in \mathcal{D}_k$, the function $\psi_k(x):=-2\phi_k
 +k(|x|^2-1)\phi_k(x) $ is an eigenfunction corresponding to $\xi_k$.}
\vs
We have an isoperimetric upper bound  for the first eigenvalue of the problem (1.19).
\vs
{\bf Theorem 1.6.} {\it Let $\Om$ be a bounded domain with smooth boundary in $\mathbf{R}^n$. Then the first nonzero eigenvalue of
the problem (1.19) satisfies
\be
\xi_1\leq \fr{(n+2)|\pa\Om|}{n|\Om|\left(\fr{|\Om|}{\omega_n}\right)^{2/n}},
\en
with equality holding if and only if $\Om$ is a ball.}
\vs
The method in proving Theorem 1.5  can be used to study some other kind of Steklov problems.
Let us consider for example the following Steklov problem
\be
\left\{\begin{array}{l}
 \ov{\Delta}^2 u= 0  \ \ {\rm in \ \ } \Om, \\
\pa_{\nu}u = \pa_{\nu}(\ov{\D} u)- \zeta \D u =0 \ \ {\rm on \ \ } \pa \Om.
\end{array}\right.
\en
The eigenvalues of this problem can be arranged as
\be
0=\zeta_0<\zeta_1\leq\zeta_2\leq\cdots\ri\infty.
\en
The smallest eigenvalue of (1.23) is zero with constant eigenfunction.  The  Rayleigh-Ritz
characterization for the first nonzero eigenvalue of (1.23) can be written as
\be
\zeta_1=\underset{\pa_{\nu}u|_{\pa\Om}=0, \int_{\pa\Om}u=0}{\underset{u\in H^2(\Om), u\neq const.}{\min}}\frac{\int_{\Om} (\ov{\D} u)^2}
{\int_{\pa \Om} |\na u|^2}.
\en
When $\Om$ is a ball, we can explicitly describe the eigenvalues and the corresponding eigenfunctions to the problem (1.23). That is, we have
\vs
{\bf Theorem 1.7.} {\it If $n\geq 2$ and $\Om =B$, then  we have

 i) the eigenvalues of (1.23) are $\zeta_0=0, \zeta_k=\fr{2k^2+nk}{k+n-2}, k=1, 2,\cdots $;

 ii) for all $k=0, 1, 2,\cdots,$ the multiplicity of $\zeta_k$ is $\mu_k$;

 iii) for all $k=0, 1, 2,\cdots,$ and  all $\phi_k\in \mathcal{D}_k$, the function $\psi_k(x):= -2\phi_k
 +k(|x|^2-1)\phi_k(x) $ is an eigenfunction corresponding to $\zeta_k$.}
\vs
The following isoperimetric upper bound  for $\zeta_1$ is an immediate consequence of (1.21), (1.25) and Theorem 1.6.
\vs
{\bf Theorem 1.8.} {\it Let $\Om$ be a bounded domain with smooth boundary in $\mathbf{R}^n$ and  $\eta_1$ the first nonzero eigenvalue of the Laplcain acting on functions on $\pa\Om$. Then the first nonzero eigenvalue of
the problem (1.23) satisfies
\be
\zeta_1\leq \fr 1{\eta_1}\cdot\fr{(n+2)|\pa\Om|}{n|\Om|\left(\fr{|\Om|}{\omega_n}\right)^{2/n}},
\en
with equality holding if and only if $\Om$ is a ball.}
\vs
An important issue in geometric analysis is to give good estimates for the  eigenvalues of various eigenvalue problems.
We refer to [4-12, 14-21, 23, 26, 27, 29] and the references therein for some interesting results about the Steklov eigenvalues.
\section{Proofs of Theorems 1.1-1.4 }
\setcounter{equation}{0}
Let us first  fix some notation. Let $M$ be $n$-dimensional compact manifold $M$ with boundary $\pa M$. We write $\la, \ra$ the Riemannian metric on $M$ as well as that induced on $\pa M$. Let $\ov{\na}$ and $\ov{\D} $ be the connection  and the Laplacian on $M$,
respectively. Let $\nu$ be the unit outward normal vector of $\pa M$. The shape operator of $\pa M$ is given by $A(X)=\na_X \nu$ and the second fundamental form of $\pa M$ is defined as $\sigma(X, Y)=\la A(X), Y\ra$, here $X, Y\in T \pa M$. The eigenvalues of $A$ are called the principal curvatures of $\pa M$ and the mean curvature $H$ of $\pa M$ is given by $H=\fr 1{n-1} {\rm trace\ } A$.
For a smooth function $f$ defined on an $n$-dimensional compact manifold $M$ with boundary $\pa M$, the following identity holds if
$h=\pa_{\nu}f|_{\pa M}$, $z=f|_{\pa M}$ and ${\rm Ric}$ denotes the Ricci tencor of $M$ (see [22], p. 46):
\be
& & \int_M \left((\ov{\D f})^2-|\ov{\na}^2 f|^2-{\rm Ric}(\ov{\na} f, \ov{\na} f)\right)\\ \no
&=& \int_{\pa M}\left( ((n-1)Hh+2\D z)h + \sigma(\na z, \na z)\right).
\en
Here $\D $ and $\na $  represent the Laplacian and the gradient on $\pa M$ with respect
to the induced metric on $\pa M$,  respectively.

We shall need the following result.
\vs
{\bf Lemma 2.1}([28]). {\it Let $(M, \la, \ra)$ be an $n$-dimensional compact connected Riemannian manifold with non-negative Ricci curvature and boundary $\pa M$.  Assume that the principal curvatures of $\pa M$ are bounded below by a positive constant $c$.  Then  the first non-zero eigenvalue  $\eta_1$ of the Laplacian  acting on functions on $\pa M$ satisfies  $\eta_1\geq (n-1)c^2$ with equality holding if and only if $M$ is isometric to an
$n$-dimensional Euclidean ball of radius $\frac 1c$.}
\vs
{\it Proof of Theorem 1.1.}
Let $u$ be an eigenfunction corresponding to $\laa_{1, \beta}$:
 Let $z$ be an eigenfunction of $\pa M$ corresponding to $\eta_1$, that is, $\D z+\eta_1 z=0$. Let  $f$ be the solution of the following  equation
\be
 \left\{\begin{array}{l}
  \ov{\D} f=0 \ \ \ \ \ \mbox{in}\ \ \  M, \\
  f|_{\pa M}=z.
\end{array}\right.
\en
Since $\int_{\pa M} z=0$, we have from (1.5)  that
\be
\laa_{1, \beta}&\leq&\fr{\int_{M} |\ov{\nabla} f|^2+\beta \int_{\pa M}|\na z|^2}{\int_{\pa M} z^2}= \eta_1\beta +\fr{\int_{M} |\ov{\nabla} f|^2}{\int_{\pa M} z^2}
\en
Setting $h=\pa_{\nu}f|_{\pa M}$, we know from (2.2) and  the divergence theorem that
\be
\int_{M} |\ov{\nabla} f|^2 =\int_{\pa M} zh.
\en
Sine the Ricci curvature of $M$ is bounded below by $-\kappa_0$ and the principal curvature of $\pa M$ are bounded below by $c$, we have
\be
{\rm Ric}(\na z, \na z)\geq -\kappa_0 |\na z|^2, \ \sigma(\na z, \na z)\geq c|\na z|^2, \ \ H\geq c.
\en
It then follows from Reilly's formula  that
\be\no
\kappa_0\int_{\pa M} zh&=&\kappa_0\int_{M}|\ov{\na} f|^2\\ \no &\geq& \int_M \left((\ov{\D} f)^2-|\ov{\na}^2 f|^2-{\rm Ric}(\ov{\na} f, \ov{\na} f)\right)\\ \no
&=&\int_{\pa M} (2(\D z)h+(n-1)h^2+\sigma(\na z, \na z)\\ \no &\geq & (n-1)c\int_{\pa M} h^2 -2\eta_1 \int_{\pa M} hz +c \int_{\pa M}|\na z|^2\\
&= & (n-1)c\int_{\pa M} h^2 -2\eta_1 \int_{\pa M} hz +c\eta_1 \int_{\pa M} z^2,
\en
that is,
\be
0\geq (n-1)c\int_{\pa M} h^2 -(2\eta_1+\kappa_0) \int_{\pa M} hz +c\eta_1 \int_{\pa M} z^2.
\en
Thus, we have
\be\no
0 &\geq& (n-1)c\int_{\pa M} \left(h -\fr{(2\eta_1+\kappa_0)}{2(n-1)c}z\right)^2 +\left(c\laa_1-\fr{(2\eta_1+\kappa_0)^2}{4(n-1)c}\right) \int_{\pa M} z^2
\\ \no &\geq& \left(c\laa_1-\fr{(2\eta_1+\kappa_0)^2}{4(n-1)c}\right) \int_{\pa M} z^2,
\en
and so
\be\no
(2\eta_1+\kappa_0)^2\geq 4(n-1)\eta_1 c^2
\en
We can also get from (2.7) that
\be
0&\geq & (n-1)c\int_{\pa M} h^2 -(2\eta_1+\kappa_0) \left(\int_{\pa M} h^2\right)^{\fr 12}\left(\int_{\pa M} z^2\right)^{\fr 12}+c\eta_1 \int_{\pa M} z^2,
\en
which, implies that
\be\no
\left(\int_{\pa M} h^2\right)^{\fr 12}\leq \frac{2\eta_1+\kappa_0+\sqrt{(2\eta_1+\kappa_0)^2- 4(n-1)\eta_1 c^2}}{2(n-1)c}\left(\int_{\pa M} z^2\right)^{\fr 12}.
\en
Combining (2.3), (2.4) and (2.9), we obtain
\be\no
\laa_{1, \beta}&\leq&  \eta_1\beta +\fr{\int_{\pa M}zh }{\int_{\pa M} z^2}\\ \no &\leq&
\eta_1\beta +\fr{\left(\int_{\pa M} h^2\right)^{\fr 12}}{\left(\int_{\pa M} z^2\right)^{\fr 12}}\\  &\leq&
\eta_1\beta +\frac{2\eta_1+\kappa_0+\sqrt{(2\eta_1+\kappa_0)^2- 4(n-1)\eta_1 c^2}}{2(n-1)c}.
\en
This proves (1.9).

If (1.9) takes  equality sign, then we have
\be\no
\left(\int_{\pa M} h^2\right)^{\fr 12}= \frac{2\eta_1+\kappa_0+\sqrt{(2\eta_1+\kappa_0)^2- 4(n-1)\eta_1 c^2}}{2(n-1)c}\left(\int_{\pa M} z^2\right)^{\fr 12}
\en
and so the inequalities in (2.6) and (2.8) should take equality sign. We infer therefore
\be
\na^2 f=0, \ \  \ \ H=c
\en
and
\be
h=\frac{2\eta_1+\kappa_0+\sqrt{(2\eta_1+\kappa_0)^2- 4(n-1)\eta_1 c^2}}{2(n-1)c} z.
\en
Take a local orthonormal fields $\{e_i\}_{i=1}^{n-1}$ tangent to $\pa M$. We conclude from (2.10) and (2.11) that
\be
0&=&\sum_{i=1}^{n-1}\ov{\na}^2 f(e_i, e_i)=\D z +(n-1)H h
\\ \no
&=& -\eta_1 z+(n-1)c\cdot\frac{2\eta_1+\kappa_0+\sqrt{(2\eta_1+\kappa_0)^2- 4(n-1)\eta_1 c^2}}{2(n-1)c}z,
\en
which gives $\kappa_0=0$ and $\eta_1=(n-1)c^2$. It then follows from Lemma 2.1 that $M $ is isometric to an $n$-dimensional Euclidean ball of radius
$\fr 1c$. On the other hand, we know from (1.7) that for the $n$-dimensional Euclidean ball of radius
$\fr 1c$, the equality  holds in (1.9). This completes the proof of Theorem 1.1.
\vskip0.3cm
{\it Proof of Theorem 1.2.} Let $ x_1,\cdots, x_n$ be the coordinate functions
on ${\mathbf R}^^n$. By choosing the coordinates origin properly, we can assume that
\be
\int_{\pa\Om} x_i =0, \ i=1,\cdots,n.
\en
It then follows from the variational characterization (1.5) for $\laa_{1, \beta}$ that for each fixed $i\in\{1,\cdots, n\}$
\be
\laa_{1, \beta}\int_{\pa \Om} x_i^2&\leq& \int_{\Om} |\ov{\na}x_i|^2+ \beta\int_{\pa \Om} |\na x_i|^2\\ \no
&=&|\Om| + \beta\int_{\pa \Om} |\na x_i|^2
\en
with equality holding if and only if $\beta\Delta x_i +\pa_{\nu} x_i = -\laa_{1, \beta} x_i$ on $\pa\Om$.

Summing over $i$ from $1$ to $n$, we get
\be
\laa_{1, \beta}\int_{\pa \Om}\sum_{i=1}^n x_i^2&\leq& n |\Om| + \beta\int_{\pa \Om}\sum_{i=1}^n |\na x_i|^2
\\ \no &=&n |\Om| + (n-1)\beta|\pa\Om|,
\en
with equality holding if and only if
\be
\beta\Delta x_i +\pa_{\nu} x_i = -\laa_{1, \beta} x_i, \ \ \ {\rm on} \ \ \pa\Om, \ \ \ \forall i\in\{1,\cdots, n\}.
\en
Take a ball $B(R,o)$ in ${\mathbf R}^n$ of radius $R$ centered at the origin so that $|B(R,o)|= |\Om|$; then
$$
R =\left(\fr{|\Om|}{\omega_n}\right)^{1/n}.
$$
By using the weighted isoperimetric inequality  in [3], we have
\be
\int_{\pa \Om}\sum_{i=1}^n x_i^2&\geq& \int_{\pa B}\sum_{i=1}^n x_i^2\\ \no
&=& |\pa B| R^2\\ \no
&=& n|\Om| \left(\fr{|\Om|}{\omega_n}\right)^{1/n}.
\en
Substituting (2.17) into (2.15), one gets (1.11). If the equality holds in (1.11), then the inequalities (2.14)
and (2.15) must take equality sign. It follows that the position vector $x = (x_1,\cdots, x_n)$ when restricted
on $\pa\Om$ satisfies
\be
\Delta x&:=&(\Delta x_1,\cdots, \Delta x_n)\\ \no &=&-\fr 1{\beta}(\pa_{\nu}x_1,\cdots,\pa_{\nu}x_n)-\fr{\laa_{1, \beta}}{\beta}(x_1,\cdots, x_n)
\\ \no &=&-\fr{1}{\beta}\nu-\fr{\laa_{1, \beta}}{\beta}(x_1,\cdots, x_n).
\en
On the other hand,  it is well known that
\be
\Delta x =(n-1){\mathbf H},
\en
where ${\mathbf H}$ is the mean curvature vector of $\pa \Om$ in ${\mathbf R}^n$. Combining (2.18) and (2.19), we have
\be
x=-\fr 1{\laa_{1,\beta}}\nu-\fr{(n-1)\beta}{\laa_{1,\beta}}{\mathbf H}, \ \ \ \ {\rm on\ \ \ }\pa \Om.
\en
Consider the function $g = |x|^2 : M\rightarrow {\mathbf R}$. It is easy to see from (2.20) that
\be\no
Z g =2\langle Z, x\rangle = 0, \ \ \ \forall Z\in {\mathfrak{X}}(\pa\Om).
\en
Thus $g$ is a constant function and so $\pa \Om$ is a hypersphere. Theorem 1.2 follows.
\vs
{\it Proof of Theorem 1.3.} i)
Let $u$ be an eigenfunction corresponding to the first eigenvalue $p_1$ of the Steklov problem:
\be
 \left\{\begin{array}{l}
  \ov{\D} u=0 \ \ \ \ \ \mbox{in}\ \ \  M, \\
  \pa_{\nu} u= p_1 u \ \ \ \ \mbox{on}\ \ \ \pa M.
\end{array}\right.
\en
Setting $w=u|_{\pa M}, y =\pa_{\nu} u|_{\pa M}$, we obtain by substituting $u$ into Reilly's formula that
\be\no
\kappa\int_ M|\ov{\na} u|^2&\geq & -2\int_{\pa M} \na w \na y + (n-1)c\int_{\pa M} y^2 + c\int_{\pa M} |\na w|^2
\\
&>& -2p_1\int_{\pa M}|\na w|^2 +c\int_{\pa M} |\na w|^2.
\en
Since $\int_{\pa M} w=0$ and $w\neq 0$, we know from the Poincar\'e inequality that
\be\no
\int_{\pa M} w^2\leq \fr 1{\eta_1}\int_{\pa M} |\na w|^2.
\en
Thus we have from the divergence theorem that
\be
\int_{M} |\ov{\na} u|^2 =\int_M {\rm div}( u\ov{\na} u)=\int_{\pa M} w y=p_1\int_{\pa M} w^2\leq \fr{p_1}{\eta_1}\int_{\pa M} |\na w|^2.
\en
Combining (2.22) with (2.23), we obtain (1.12).

 ii) Let $f$ be an eigenfunction corresponding to the first eigenvalue $\laa_{1, \beta}$ of the Wentzell-Laplace operator
of $M$. Setting
$$
\gamma=\fr 1{\beta},\ \laa=\fr{\laa_{1, \beta}}{\beta},\ z=f|_{\pa M},\ h=\pa_{\nu} f|_{\pa M}, $$
we have
\be
\ov{\Delta} f =0
\en
and
\be
 \Delta z =\gamma h-\laa  z.
\en
It then follows from  Reilly's formula that
\be\no
0&\geq& \int_M \left((\ov{\D} f)^2-|\ov{\na}^2 f|^2-{\rm Ric}(\ov{\na} f, \ov{\na} f)\right)\\ \no
&=&\int_{\pa M}(2(\D z)h +(n-1)Hh^2+\sigma(\na z, \na z))\\ \no &\geq &\int_{\pa M}(2(\gamma h-\laa z)h+ (n-1)ch^2 +c |\na z|^2)\\ \no
&= & \int_{\pa M}(2(\gamma h-\laa z)h+ (n-1)ch^2 -c z\D z)\\ \no
&= & \int_{\pa M}( c\laa z^2-(2\laa+c\gamma)zh+(2\gamma+(n-1)c) h^2)\\
&=& (2\gamma+(n-1)c)\int_{\pa M}\left(h-\fr{2\laa+c\gamma}{2(2\gamma+(n-1)c)}z\right)^2\\ \no & & +\left(c\laa -\fr{(2\laa+c\gamma)^2}{4(2\gamma+(n-1)c)}\right)\int_{\pa M} z^2
\\ \no
&\geq &\left(c\laa -\fr{(2\laa+c\gamma)^2}{4(2\gamma+(n-1)c)}\right)\int_{\pa M} z^2,
\en
which, gives
\be\no
c\laa -\fr{(2\laa+c\gamma)^2}{4(2\gamma+(n-1)c)}\leq 0.
\en
Thus, we have either
\be
\laa\geq \fr{\left(\gamma+(n-1)c+\sqrt{(n-1)^2c^2+2(n-1)\gamma c}\right)c}{2},
\en
or
\be
\laa\leq \fr{\left(\gamma+(n-1)c-\sqrt{(n-1)^2c^2+2(n-1)\gamma c}\right)c}{2}.
\en
We {\it claim} that (2.28) does not occur. In fact, multiplying (2.25) by $-z$ and integrating on $\pa M$, we get
by using the divergence theorem and (2.24) that
\be
\laa \int_{\pa M} z^2& =&\int_{\pa M} |\na z|^2 +\gamma \int_{\pa M} h z
\\ \no & =&\int_{\pa M} |\na z|^2 +\gamma \int_{M} |\ov{\na}f|^2.
\en
Sine $\int_{\pa M} z=0, z\neq 0$, we have from the Poincar\'e inequality and the variational characterization (1.5) that
\be
\int_{\pa M} |\na z|^2 \geq \eta_1 \int_{\pa M} z^2,\  \int_{M} |\ov{\na}f|^2\geq p_1 \int_{\pa M} z^2.
\en
Also, we know from (1.12) that
\be
p_1 > \fr c2.
\en
From (2.29)-(2.31) and  Lemma 2.1, we conclude that
\be
\laa > \eta_1 +\fr{c\gamma}2\geq (n-1)c^2+\fr{c\gamma}2,
\en
which proves our {\it claim}. In order to finish the proof of Theorem 1.2, we need only to exclude the equality case in (2.27).
We shall do this by contradiction. Thus suppose that
\be
\laa= \fr{\left(\gamma+(n-1)c+\sqrt{(n-1)^2c^2+2(n-1)\gamma c}\right)c}{2}.
\en
From the proof above, we know that all the
inequalities in (2.26) should take equality sign. It then follows that the principal curvatures of $\pa M$ satisfy
\be
\rho_1=\cdots=\rho_{n-1}=c \ \ {\rm on \ \ } M,
\en
that
\be
\ov{\na}^2 f=0 \ \ {\rm on \ \ } M
\en
and that
\be
h=\fr{2\laa+c\gamma}{2(2\gamma+(n-1)c)}z  \ \ {\rm on \ \ } \pa M.
\en
Restricting (2.35) on $\pa M$ and using (2.34), we conclude that
\be
 \ h=cz.
\en
Combining (2.36) with (2.37), we know that
\be
\fr{2\laa+c\gamma}{2(2\gamma+(n-1)c)}=c.
\en
We infer  from (2.33) and (2.38) that $\gamma=0$. This is a contradiction. The proof of Theorem 1.2 is complete.
\vs
Theorem 1.4 is a special case of a more general result. Before stating it, let's introduce the following
\vskip0.2cm
{\bf Definition 2.1.} {\it Let $M$ be $n$-dimensional compact manifold $M$ with boundary and $\tau$ a fixed positive number. We define the
$\tau$-fundamental tone of $M$ to be
\be
F_{\tau, M}=\underset{\int_{\pa \Om}u=0}{\underset{0\neq u\in H^2(\Om)}{\min}}\frac{\int_{\Om}( |\ov{\nabla}^2 u|^2+ \tau |\ov{\nabla} u|^2)}
{\int_{\pa \Om} u^2}.
\en
}
\vskip0.2cm
We have the following result  which implies Theorem 1.4.
\vskip0.2cm
{\bf Theorem 2.2.} {\it Let $M$ be an $n$-dimensional compact manifold with  boundary and assume that the principle curvatures of $\pa M$ are bounded below by a positive constant $c$. Let $\eta_1$ be the first  non-zero eigenvalue of the Laplacian  acting on functions on $\pa M$.
If the Ricci curvature of $M$ is bounded below by $-\kappa$ for some constant $\kappa\geq 0$, then the $\tau$-fundamental tone of $M$  satisfies
\be
F_{\tau, M}\leq\frac{(2\eta_1+\kappa)^2+ 2\tau\left(2\eta_1+\kappa+\sqrt{(2\eta_1+\kappa)^2- 4(n-1)\eta_1 c^2}\right)}{4(n-1)c}-c\eta_1,
\en
with equality holding if and only if $\kappa=0$ and $M$ is isometric to a ball of radius $\frac 1c$ in ${\bf R}^n$.
}
\vskip0.2cm
{\it Proof of Theorem 2.2.}  As in the proof of Theorem 1.1, let $z$ be an eigenfunction of $\pa M$ corresponding to $\eta_1$ and $f$  the solution of the  equation
\be
 \left\{\begin{array}{l}
  \ov{\D} f=0 \ \ \ \ \ \mbox{in}\ \ \  M, \\
  f|_{\pa M}=z.
\end{array}\right.
\en
 Setting $h=\pa_{\nu}f|_{\pa M}$, we have from the definition of $F_{\tau, M}$  that
\be\no
F(\tau, M)&\leq& \fr{\int_{M} |\ov{\nabla}^2 f|^2+ \tau\int_{M}|\ov{\na} f|^2}{\int_{\pa M} z^2}\\
&=&\fr{\int_{M} |\ov{\nabla}^2 f|^2+ \tau\int_{M}zh}{\int_{\pa M} z^2}.
 \en
Since
\be
{\rm Ric}(\ov{\na}f, \ov{\na}f)\geq -\kappa|\ov{\na}f|^2,\  H\geq c,\  \sigma(\na z, \na z)\geq c|\na z|^2,
\en
we know  from Reilly's formula that
\be\no
-\int_M |\ov{\na}^2 f|^2+\int_M \kappa|\ov{\na}f|^2
&\geq &(n-1)c\int_{\pa M} h^2 -2\eta_1 \int_{\pa M} hz +c \int_{\pa M}|\ov{\na} z|^2\\ \no
&= & (n-1)c\int_{\pa M} h^2 -2\eta_1 \int_{\pa M} hz -c\eta_1 \int_{\pa M} z^2,
\en
which gives
\be\no
\fr{\int_M |\ov{\na}^2 f|^2}{\int_{\pa M} z^2}&\leq& -(n-1)c\fr{\int_{\pa M} h^2}{\int_{\pa M} z^2}+(2\eta_1+\kappa)\left(\fr{\int_{\pa M} h^2}{\int_{\pa M} z^2}\right)^{1/2}-c\eta_1\\
&\leq&-c\eta_1+\fr{(2\eta_1+\kappa)^2}{4(n-1)c}.
\en
Also, as in the proof of Theorem 1.1, we have (Cf. (2.9)), we have
\be
\fr{\int_{M}zh}{\int_{\pa M} z^2}\leq \frac{2\eta_1+\kappa+\sqrt{(2\eta_1+\kappa)^2- 4(n-1)\eta_1 c^2}}{2(n-1)c}.
\en
Substituting (2.44) and (2.45) into (2.42), one gets (2.40).

If the equality  holds in (2.40), we can use the same arguments as in the proof of Theorem 1.1 to conclude that $\eta_1=(n-1)c^2$. Thus, $M$ is isometric to a ball of radius $\fr 1c$ by Lemma 2.1. Also, if $M$ is a ball of radius $\fr 1c$ in ${\mathbf R}^n$, we have $\eta_1=(n-1)c^2$,
${\rm Ric}_M=0, F_{\tau, M}= \fr{\tau}c$ and therefore, the equality holds in (2.40). This completes the proof of Theorem 2.2.
\section{Proofs of Theorems 1.5-1.8}
\setcounter{equation}{0}
In this section, we shall prove Theorems 1.5-1.7.
\vs
{\it Proof of Theorem 1.5.} Let $u$ be an eigenfunction of (1.19) with $\Om=B$ corresponding to an eigenvalue $\xi$. Since $u$ is bi-harmonic, there  exist uniquely determined harmonic functions
$g, h$ on $B$, such that $u(x)= g(x)+|x|^2 h(x)$ (Cf. [1]). The Euler operator $\La$ will play a key role in the proof :
$$\La f(x):=\sum_{i=1}^n x_i\fr{\pa f}{\pa x_i}(x).$$
We will use the notation $\La^lf=\La(\La^{l-1}f), l=1, 2,\cdots. $
Observe that
 $\La$ is related to the exterior normal derivative, $\pa_{\nu}f|_{\pa B}=\La f|_{\pa B}$, and
if $f$ is harmonic then so is $\La f$. Setting $g+h= -2z$, we have
$u(x)=-2z(x)+(|x|^2-1)h(x)$. From
$$0=\pa_{\nu} u|_{\pa B}=(-2\pa_{\nu} z+ 2h)|_{\pa B}=(-2\La z+ 2h)|_{\pa B}$$
and the fact $-2\La z+ 2h$ is harmonic, we know that $-2\La z+ 2h=0$ on $B$ and so we have
$$u(x)=-2z(x)+(|x|^2-1)\La z(x).$$
It is easy to see that
\be\no
\ov{\D} u= 2n\La z +4\La^2z.
\en
Thus
\be
\pa_{\nu}(\ov{\D} u)|_{\pa B}=(2n \La^2z +4\La^3z)|_{\pa B}=-\xi u|_{\pa B}=2\xi z|_{\pa B},
\en
which gives the following important relation:
\be
n \La^2z +2\La^3z= \xi z \ \ {\rm on \ } B.
\en
Since $z$ is harmonic on $B$, there exist $p_m\in {\mathcal D}_m$ such that
\be
z(x)=\sum_{m=0}^{\infty} p_m(x)
\en
for all $x\in B$, the series converging absolutely and uniformly on compact subsets of $B$(Cf. [2]). Substituting (3.3) into (3.2) and observing
$\La p_m = mp_m$, we get
\be
2m^3+ nm^2 = \xi, m=0, 1,\cdots.
\en
This implies that  all but one of the $p_{m}'s$ are zeros. Thus the eigenvalues  are $\xi_k=2 k^3+n k^2, k=0, 1, \cdots, $ and the multiplicity of $\xi_k$
is the dimension of ${\mathcal D}_k$. Also, we know from the above proof that for all $\phi_k\in {\mathcal D}_k$, the function $\psi_k(x) :=-2\phi_k
 +k(|x|^2-1)\phi_k(x)$ is an eigenfunction corresponding to $\xi_k$. This completes the proof of Theorem 1.5.
\vskip0.3cm
{\it Proof of Theorem 1.6.} Let $x_i, i=1,\cdots,n,$ be the coordinate functions on ${\mathbf R}^n$. By taking the coordinates origin properly, we can assume that
\be
\int_{\Om} x_i=0, \ i=1,\cdots,n.
\en
For each $i\in \{1,\cdots, n\}$, let $g_i$ be the solution of the problem
\be
 \left\{\begin{array}{l}
  \ov{\D} g_i= x_i\ \ \ \ \ \mbox{in}\ \ \  \Om, \\
  \pa_{\nu}g_i|_{\pa \Om}=0.
\end{array}\right.
\en
We can assume without lose of generality that
\be
\int_{\pa \Om} g_i=0, \ \ i=1,\cdots,n.
\en
By using the Rayleigh-Ritz characterization (1.21), we get
\be
\xi_1\leq \fr{\int_{\Om} x_i^2}{\int_{\pa \Om}  g_i^2}, \ \ i=1,\cdots,n.
\en
By using (3.6) and the divergence theorem we get
\be
\int_{\Om} x_i^2 =\int_{\Om} x_i \ov{\D} g_i=-\int_{\Om}\la \ov{\na} x_i, \ov{\na} g_i\ra=-\int_{\pa \Om} g_i\pa_{\nu}x_i,
\en
which implies that
\be
\left(\int_{\Om} x_i^2\right)^2\leq\int_{\pa\Om} (\pa_{\nu} x_i)^2\int_{\pa\Om} g_i^2.
\en
Substituting (3.10) into  (3.8), we infer
\be
\xi_1\leq  \fr{\int_{\pa \Om}(\pa_{\nu} x_i)^2 }{\int_{\Om} x_i^2}, i=1,\cdots,n.
\en
Summing over $i$, we have
\be
\xi_1 \int_{\Om} \rho^2\leq \int_{\pa\Om} \sum_{i=1}^n (\pa_{\nu} x_i)^2 =|\pa \Om|,
\en
where $\rho= d(o,\cdot): \Om\ri {\mathbf R}$ is the distance function from the origin. Let $B(R, o)$ be the ball of radius $R$  centered at the origin in ${\mathbf R}^n$
such that $|B(R, o)|=|\Om|$ and set $W=\Om\cap B(R, o)$; then
\be\no
\int_{\Om} \rho^2 &=&\int_{W} \rho^2 + \int_{\Om\setminus W} \rho^2\\ \no
&\geq&\int_{W} \rho^2 + \int_{\Om\setminus W} R^2\\ \no
&=& \int_{W} \rho^2 + R^2|B(R, o)\setminus W|\\
 &\geq &\int_W \rho^2 + \int_{B(R, o)\setminus W} \rho^2\\ \no
 &=& \int_{B(R, o)} \rho^2\\ \no
 &=& \fr{n\omega_{n}}{n+2} R^{n+2}\\ \no
 &=& \fr{n\omega_{n}}{n+2}\left( \fr{|\Om|}{\omega_n}\right)^{\fr{n+2}n},
 \en
 which, combining with (3.12), gives (1.22). If the equality holds in (1.22), then the inequalities (3.10)-(3.13)  should take equality sign, which easily implies  that
 $\Om$ is a ball. On the other hand, we know by scaling and Theorem 1.5 that
 \be\xi_1(B(R, o))=\fr{n+2}{R^3} =\fr{(n+2)|\pa B(R, o)|}{n|B(R, o)|\left(\fr{|B(R, o)|}{\omega_n}\right)^{2/n}},\en
that is, the equality  holds for balls in (1.22). This completes the proof of Theorem 1.6.
\vs
{\it Proof of Theorem 1.7.} Let $f$ be an eigenfunction of (1.23) corresponding to a nonzero eigenvalue $\zeta$. Since $f$ is bi-harmonic with $\pa_{\nu}f|_{\pa B}=0$, by using the same arguments as in the proof of Theorem 1.5, there exists a  uniquely determined harmonic function
$w$ on $B$, such that
$$f(x)= -2w(x)+(|x|^2-1) \La w(x).$$
We have
\be
\ov{\D} f=2n \La w + 4\La^2w
\en
and if $y\in \pa B$,
\be\no
\ov{\D} f(y)&=&\D f(y) +\ov{\na}^2(\nu, \nu)(y)\\ \no &=&\D f(y) + (\La^2f-\La f)(y)\\
&=& \D f(y) + \La^2f(y).
\en
Also, we have
\be
\pa_{\nu}(\ov{\D} f)|_{\pa B}=(2n \La^2 w +4\La^3w)|_{\pa B}
\en
It is easy to see that for any $x\in B$,
\be
\La f(x)&=&-2\La w(x)+\sum_{i=1}^n x_i((|x|^2-1)\La w)_{x_i}\\ \no &=&(|x|^2-1)(2\La w +\La^2 w))(x)
\en
and
\be \no
\La^2f(x)= 2|x|^2 (2\La w +\La^2 w)(x)+(|x|^2-1)(2\La^2 w +\La^3 w)(x).
\en
Hence
\be
\La^2f(y)= 2 (2\La w +\La^2 w)(y), \ \ \forall y\in \pa B,
\en
which, combining with (3.15) and (3.16), gives
\be
\D f(y)=((2n-4)\La w+2\La^2 w)(y), \ \ \forall y\in \pa B.
\en
It then follows from $(\pa_{\nu}(\ov{\D} f)-\zeta \D f)|_{\pa B}=0$, (3.17) and (3.20) that
\be
(2\La^3 w+n\La^2w-\zeta(\La^2 w+(n-2)\La w))|_{\pa B}=0.
\en
Consequently, we have
\be
2\La^3 w+n\La^2w=\zeta(\La^2 w+(n-2)\La w)\ \ \ {\rm on} \ \ B.
\en
Using the same arguments as in  the final part of the proof of Theorem 1.5, we can conclude the conclusions of Theorem 1.7.
\vs
{\it Proof of Theorem 1.8.} Let $u$ be an eigenfunction corresponding to the first eigenvalue $\xi_1$ of the problem:
\be
\left\{\begin{array}{l}
 \ov{\Delta}^2 u= 0  \ \ {\rm in \ \ } \Om, \\
\pa_{\nu}u = \pa_{\nu}(\ov{\D} u)+\xi_1 u =0 \ \ {\rm on \ \ } \pa \Om;
\end{array}\right.
\en
then
\be
\xi_1=\fr{\int_{\Om}(\ov{\D} u)^2}{\int_{\pa \Om}u^2}.
\en
It follows from the variational characterization (1.25) for $\zeta_1$ that
\be
\zeta_1\leq \fr{\int_{\Om}(\ov{\D} u)^2}{\int_{\pa \Om}|\na u|^2}.
\en
We have from the Poincar\'e inequality that
\be
\int_{\pa \Om}|\na u|^2\geq \eta_1\int_{\pa \Om}u^2,
\en
with equality holding if and only if $u|_{\pa \Om}$ is an eigenfunction corresponding to the eigenvalue $\eta_1$. Combining (3.24)-(3.26) with Theorem 1.6, one can finish the proof of Theorem 1.8.

 \vskip0.4cm
\noindent Changyu Xia ( xia@mat.unb.br )

\noindent Departamento de Matem\'atica, Universidade de Bras\'{\i}lia, 70910-900 Bras\'{\i}lia-DF, Brazil

\vs
\noindent Qiaoling Wang ( wang@mat.unb.br )

\noindent Departamento de Matem\'atica, Universidade de Bras\'{\i}lia, 70910-900 Bras\'{\i}lia-DF, Brazil

\end{document}